# Making Mathematical Claims

*How many non-congruent quadrilaterals can you make with vertices in the 9-dot square grid shown in Figure 1? How would you organize those quadrilaterals into categories?[1]*

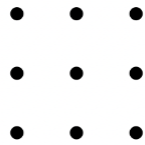

Figure 1

*What are all the ways to reflect the pirate face along lines $\ell_1, \ell_2, \ell_3, \ell_4$ so the pirate lands in the ocean (the shaded region in Figure 2)? Which ways would you classify as equivalent? Inequivalent? Why?*

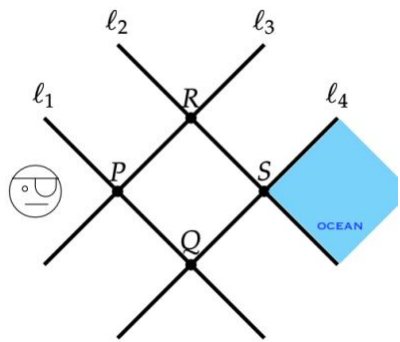

Figure 2

Good problems grab us. They invite us to find patterns, make conjectures, and prove — or perhaps disprove — a conjecture. When I first taught, I saw my work as tantalizing students with structures just beyond their reach, so that I could elicit conjectures from promising half-phrases. With a community conjecture crystallized on the board, "we"[2] proved the statement. "We" anointed the conjecture a community theorem, and "we" moved on. I hoped that, through repeated exposure to this routine, students would absorb a mathematical process from discovery to proof.

But I've since wondered: What does this routine teach students?

I've concluded that if this is the only instructional routine that students experience, they may leave with an impoverished image of the beauty and joy that doing math can offer. Moreover, to nourish students' mathematical participation, we must find ways to cultivate their mathematical language beyond modeling precision ourselves. We must find ways to develop what the linguist M. A. K. Halliday termed a *mathematical register*: a set of meanings that are cued by mathematical use of language, and that are crucial to express for mathematical purposes. Without experience with much mathematical expression, formal mathematical writing can seem stilted. We write thickets of conjunctions and compound noun phrases with implicit logical relationships.[3] Try as we might to be simple and direct,[4] stripping away intricacies only seems to lose intended mathematical meaning. How are students to appreciate, let alone be willing to express their own thinking, in this kind of writing?

---

[1] This problem is due to Pat Janike, a high school teacher for Lincoln Public Schools. He developed this problem as an opener for a unit on congruence.

[2] In this essay, I use *"we"* when I as an instructor have effectively spoken for myself and my students, and *we* when expressing a view that I believe is commonly held in the professional mathematical community.

**Appreciating a mathematical register**

Last fall semester, I tried a brilliant activity by Sam Shah.[5] He developed it for high school geometry and I adapted it for a Modern Geometry course for prospective high school teachers. The activity, *Attacks and Counterattacks in Geometry*, comes in three parts and a sequel.
- **Part One.** Students are asked to define, individually and then in groups, the terms "triangle", "circle", and "polygon". This is the "attack".
- **Part Two.** After trading definitions, the groups are asked to "counterattack" each other's definitions, that is, to come up with an example of something that satisfies the definition they are reading but is not actually a triangle (or circle, or polygon).

After group members shared their individual definitions with each other, I asked groups to write a consensus definition on whiteboards in a (topological) circle around the room. Then, each group walked counterclockwise to the next whiteboard to suggest counterattacks and questions (see Figure 3).

| Object: | Triangle | Circle | Polygon |
|---|---|---|---|
| Attack: | A shape with three straight sides who has three angles that sum to 180° | A closed shape with 1 continuous side that returns to the same point by a curved path | A 2-dimensional shape with 3 or more straight edges and three or more angles |
| Counterattack: | 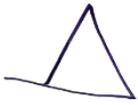 | 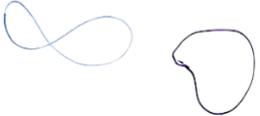 | 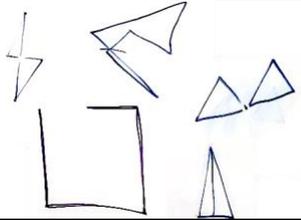 |

**Figure 3.** Prospective high school teachers' "attacks" (proposed definitions) and

---

[3] See Schleppegrell (2007) for one account of the intricacy of mathematical language as compared to natural language. Perhaps the price of the "unreasonably effectiveness of mathematics" (Wigner, 1960) are the linguistic contortions needed for precision and hierarchical logical dependencies.
[4] Indeed, Barzun's (1975) guide to writing is titled *Simple & Direct: A Rhetoric for Writers*.
[5] https://samjshah.com/2014/10/19/attacks-and-counterattacks-in-geometry/

"counterattacks" (object that satisfies proposed definition in an unintended way).

Two student groups identified proposed definitions that seemed impermeable to counterattacking. Each group called other students over to help counterattack, and soon, heated discussion was taking place in front of these two proposed definitions. Each time a student suggested a counterattack, another student pointed out why the counterattack did not meet the conditions of the proposed definition. I asked students to return to their desks, for the next part of the activity:

- **Part Three.** Students are given copies of "textbook definitions" of the terms: correct definitions that are quoted from textbooks. I provided two different and mathematically equivalent definitions of each term, and told the students that these definitions were correct.

To my astonishment, my students faced down these definitions with attempts to counterattack. After accepting that these definitions "worked", they critiqued the phrasing and suggested ways to reorder the sentences or use different phrases for clarity. Students debated the resulting meaning of the suggestions. We finished by drafting a class community definition that they voted unanimously to approve, in both meaning and clarity. (We also concluded that the two definitions they debated about were in fact correct.) Finally, the activity ended:

- **Sequel.** Students are asked to counterattack the "textbook definitions" — with certain clauses crossed out. (The students found this part relatively easy, and it was a nice way to wrap up.)

These students had found their "mathematical register".[6] They read and wrote mathematical text with meaning, and they understood that text was theirs to approve, disapprove, and improve. I continued to use the attacks/counterattacks routine, in abbreviated ways, throughout the term. For the first time in 15+ years of teaching collegiate mathematics, students were excited to receive definitions, and multiple students suggested improvements to definitional texts.

## Claim making

Willing self-expression in a mathematical register is prerequisite to agency in making formal claims. As the students tapped into the power of a mathematical register, they freely and respectfully took up and critiqued each other's conjectures. "It would be shorter if you said it this way" — suggesting a synthesis of conditions. "But if this is true, can't you also say this?" — suggesting a generalization.

Part of what made these conversations possible was the idea of "claim". Some years ago, when I was interested in making mathematical processes more explicit, I explained that:

*Mathematically, a "claim" is a statement that is provable or disprovable (or provable to be unprovable) through a deductive logical argument.*

However, I quickly gathered that while this description may suffice for an audience of those already in the know, it is not enough for teaching. Consider students doing the 9-dot grid problem that opens this essay. Under the above conception of "claim", statements such as "The array has 9 dots", or "There are 36 ways to choose pairs of points from the grid" are claims. To be sure, they aren't going to advance collective thinking on the problem, but they are claims nonetheless. And so, I propose the following conception of "claim" for the purposes of teaching:

> *A "claim" is always in reference to a question. It is a statement that answers the question in a potentially satisfying way, and that can be proven to be true or false (or proven to be unprovable).*

---

[6] A "register" can be thought of as a linguistic signature for a specialized, professional use.

Sometimes, I use the "preschooler test" to explain the meaning of "satisfying". Suppose a preschooler asks you whether all stop signs are red. If you answer a different question, the preschooler will be unhappy. It is also not satisfying to respond, "The one in front of us is red." A more satisfying answer might explain that red can signal alertness, and that in many countries, stop signs are red. Satisfying claims speak to a question, indicate why an underlying pattern makes sense, and are as general as possible. (In reality, there is probably nothing that can satisfy a curious preschooler, but I find that this conceit is helpful in teaching nonetheless.) The use of "satisfying" here taps into two senses of the word: mathematical and emotional. The claim should provide a mathematical answer to the question, and in a way that appeals to a curious question poser.

I also explain that we can think of claims as an "I bet" statement.[7] If you are the arbitrator for a bet between two friends, you would want to make absolutely sure that everyone knows exactly what the statement means, and also what evidence and inferences are acceptable for determining the truth of the statement.

**Claim proving**

Sometimes, a claim is small. For instance, I asked students to determine whether the following statement is true or false:

> *Let c be a constant nonzero real number. Plot all points (x, y) such that the ratio $\frac{y}{x} = c$. Then you have plotted a line.*

I confess that I assigned this problem thinking that the statement was true. Most students blithely agreed. But then, one group grew louder and louder. A student pounded, "No, look, I have a *proof*. This claim is false." His group mates commented on his boldness— was he serious he had a *proof*? The interaction turned students' heads in every other group. He pointed out that when *x* = 0, the ratio is undefined, and when *y* = 0, there is no solution for nonzero *c*. The points plotted formed a line missing two points. Some protested, "But that's basically a line!" Another rejoined, "But what's the definition of a line?" Everyone laughed; asking after "definitions" had become something of a friendly inside joke by this point. The class quieted and agreed: the statement was false, but would be true if the equation read *y* = *cx* rather than $\frac{y}{x} = c$. With this alternative equation, the locus of points would satisfy – mathematically – the definition of a line. These students took the notion of a "satisfying answer" to go beyond "true" or "false"; they sought to find a true statement out of a false one, and to delve into the details of what made the original statement false. There was an air of collective satisfaction – in the emotional

---

[7] I learned this formulation from Annie Selden.

sense – when they arrived at the insight about the equation *y = cx*. No one student came up with all the answers; they engaged in the exploration together, spurred by students' commitment, camaraderie, and curiosity.

Other times, a claim is larger, such as when figuring out how these six statements could all define an ellipse.

1. *Graph for an equation of the form $\frac{x^2}{a^2} + \frac{y^2}{b^2} = 1$, where $a, b \neq 0$.*
2. *A stretched or squished circle.*
3. *A shape you can draw by attaching a string to two thumbtacks and pulling the string as far as it will go with a pencil tip, and then marking all the points the pencil tip goes.*
4. *Given a real number c, it is the set of (x, y) such that the sum of the distances from (x, y) to (− c, 0) and (c, 0) is 2a, where a is a positive constant.*
5. *A slice of a hollow cylinder.*
6. *A slice of a hollow cone.*

When I first designed this activity a few years ago, bi-infinite cylinders annoyed students, as did slices parallel to a bi-infinite cylinder's axis, which are the only slices that are not ellipses. "We" made conjectures that involved my imposition of exceptional cases. As students set to proving conjectures, they rolled their eyes when I asked how their reasoning took into account exceptional cases.

This time, the class expressed skepticism that bi-infinite cylinders should be a standard concept, but bought into it when they saw it made conjectures easier to write. They generated the exceptional cases themselves, and edited their own budding conjectures accordingly. When I later asked how their proofs accounted for the non-elliptical slices, some students did sigh at the potential difficulty of the question, but they engaged with it nonetheless. Some students also asked what happens with imaginary values for *a* and *b* in the equation, and we explored hyperbolas.

Later on, when we did the pirate task opening this essay, students conjectured that the composition of two reflections is a translation, rotation, or the identity transformation. Although we needed to develop some machinery to prove this conjecture, the process went more smoothly because the conjecture came from the students. Their willingness to participate in a mathematical register continued to pay dividends.

**Reconceiving routines around mathematical claims**

A good problem that invites exploration can lead to conjecture and proof, because in theory, "try something out and see what happens" can lead to conjecture. Yet if students are unable to express a conjecture, in all its precision, then they may not be participating fully. Over the past few years, I have tried to shift students' frame for exploration from

*We are exploring*

to

*We are exploring for the purpose of claim making, which is for the purpose of claim proving, disproving, and refinement.*

In previous years, I have mostly failed at changing the frame, despite reflections on mathematical processes, stating my intentions explicitly, and giving much scaffolding for specific conjectures and associated proofs. None of these interventions seemed to make a significant difference. Moreover, I worried that scaffolding reinforced that it is the teacher's responsibility to tell students when to explore, conjecture, and prove.

As mathematicians and teachers, one of our main goals is for students to do mathematics. Part of doing mathematics is self-regulating shifts from exploration to claim making to reasoning about claims and back. My theory now is that for students to take charge of the mathematics, we must set up students to tap into a mathematical register, and define claims not for mathematics but for engaging students in doing mathematics. A mathematical register is needed to understand the precise implications of a statement, which is needed to express oneself as well as interact with others about the mathematics. Introducing the notion of a "satisfying" answer supported social accountability for advancing explorations.

I used to frame the mathematical proving process as

1. *exploring;*
2. *claim making; and*
3. *claim proving*

as three distinct phases, where students' active participation ended with claim making. In retrospect, I positioned conjectures as an epistemological end. Although I still do this at times, I now try to find ways to cultivate an alternative frame:[8]

0. *finding something out*
1. *exploring to propose a provable claim that helps with finding out about that something;*
2. *proving, disproving, refining, or generating alternate claims; and*
3. *launching new phenomena to discover.*

As I begin this year of teaching, I'm looking forward to seeing what this year's students will teach me about teaching mathematics, and how my thinking about claim making and a mathematical register for the purposes of teaching will continue to evolve.

**Acknowledgments.** I am grateful to Sam Shah for blogging his ideas for teaching, to Erin Baldinger for ongoing support in my teaching and thinking about teaching, to Allan Donsig for his openness to our department members' trying new things in teaching, and to the undergraduates who tried on new roles as students as I tried on a new role as teacher.

---

[8] This frame is a variation on Hyman Bass's (2015) description of the phases of mathematical discovery as an iterative dynamic among exploration, discovery, conjecture, proof, verification.